\documentclass[11pt,twoside]{amsart}
        \usepackage{amssymb}
\usepackage{graphicx,color}
\usepackage{cleveref}
\usepackage{todonotes}

\newtheorem{thm}{Theorem}[section]

\newtheorem{conjecture}[thm]{Conjecture}
\newtheorem{problem}[thm]{Problem}
\newtheorem{question}[thm]{Question}

\newcommand{\skipit}[1]{{}}
\newcommand{\prfend}{\hbox to7pt{\hfil}
\par\vskip-\baselineskip\hbox to\hsize
{\hfil\vbox {\hrule width6pt height6pt}}\vskip\baselineskip}

\newcommand {\PP}{\mathbb{P}}

\newcommand{\myarrow}[2]{\hbox to #1pt{\hfil$\to$\hfil}{\hskip-#1pt{\raise
10pt\hbox to#1pt{\hfil$\scriptscriptstyle #2$\hfil}}}}

\begin{document}

\title[List of problems]
{List of problems}
 \author[M. Juhnke-Kubitzke]{Martina Juhnke-Kubitzke}
 \address{Universit\"at Osnabr\"uck, Albrechtstra\ss e 28a, 49076 Osnabr\"uck, Germany}
  \email{juhnke-kubitzke@uni-osnabrueck.de}.
  \author[R.\ M.\ Mir\'o-Roig]{Rosa M.\ Mir\'o-Roig}
  \address{Facultat de
  Matem\`atiques i Inform\`atica, Universitat de Barcelona, Gran Via des les
  Corts Catalanes 585, 08007 Barcelona, Spain} \email{miro@ub.edu, ORCID 0000-0003-1375-6547}

\thanks{The second author has been partially supported by the grant PID2019-104844GB-I00}

\begin{abstract} The study of the Lefschetz properties of Artinian graded algebras was motivated by the hard
Lefschetz theorem for a smooth complex projective variety, a breakthrough in algebraic topology and geometry.    Over the last few years, this topic has
attracted increasing attention from mathematicians in various areas.
Here, we suggest some important open problems about or related to Lefschetz
properties of Artinian graded algebras  with the ultimate aim to attract the attention of young researchers from different areas.

\end{abstract}

\maketitle

\section{Introduction}

As it was pointed out in the previous chapters of this book, the study of  Lefschetz properties of Artinian algebras was originally motivated by the Lefschetz theory for projective manifolds, begun by S. Lefschetz and well-established by the late 1950s. Many of the important Artinian graded algebras appear as cohomology rings of an algebraic variety or manifold, though recent important developments have demonstrated cases of the Lefschetz property beyond such geometric settings (such as Coxeter groups or matroids). Lefschetz properties also appear as one important ingredient of the K\"ahler package.
In the last two decades there has been fascinating progress on the
study of the weak and strong Lefschetz property from different perspectives, inspired in part by, and contributing to, developments in algebraic geometry, commutative algebra and combinatorics among others
but we want to emphasize that in spite of the big progress in the area made during these last decades a lot of interesting problems remain open. Lefschetz properties have shown to be ubiquitous since this subject matter has connections to many branches of mathematics. Indeed, a central object of study are Gorenstein algebras (also known as Poincar\'e duality algebras) which are of strong interest not only in algebraic geometry but also in commutative algebra, algebraic topology and combinatorics. Notably, several important results in this area have been obtained by using unexpected methods and finding unexpected connections between apparently different topics.

In this last chapter, we gather a collection of open problems about or related to Lefschetz properties  with the ultimate aim to attract the attention of young researchers from different areas. Many researchers have contributed to this list which in no means claims to be an exhaustive list but which hopefully gives a flavor of the problems in this area. Moreover, the suggested problems vary significantly in their level of detail.

\vskip 2mm
In the following, we divide the problems into 4 blocks:
\begin{itemize}
    \item[(1)] Failure/presence of Lefschetz properties for special types of algebra.
    \item[(2)] Geometric aspects related to Lefschetz properties.
    \item[(3)] Lefschetz properties in combinatorics.
    \item[(4)] Lefschetz properties and Jordan types of Artinian algebras.
\end{itemize}.

\section{Failure/presence of Lefschetz properties for special types of algebra}

Though many algebras are expected to have the weak Lefschetz property (WLP, for short), establishing this property is often rather difficult. In this first section we gather some open problems related to the failure/presence of the WLP for special types of Artinian algebras. We also mention that some of the main results that have been achieved up to date only hold over fields of characteristic zero or under other restrictions on the characteristic of the ground field. It is hence of great interest to understand the influence of the characteristic of the ground field on Lefschetz type questions. If not otherwise stated, in this section, we will always assume to work over a field $\mathbb{K}$ of characteristic zero.

\vskip 4mm
{\bf Complete Intersections.}  It is known that  all graded algebras in two or fewer variables have the WLP \cite[Proposition 4.4]{HMNW}.  In particular, these includes height two complete intersections in two or fewer variables. Following up on our previous comment on the characteristic dependence, we want to mention that it is not hard to see that complete intersections in two variables can fail the WLP in positive characteristic. Coming back to characteristic zero, Harima, Migliore, Nagel and Watanabe \cite{HMNW}  showed that any Artinian height three complete intersection has the WLP. On the other hand, Stanley \cite{stanley}  and Watanabe  \cite{watanabe} showed that \emph{any} monomial complete intersection, independent of the number of variables, has even the strong Lefschetz property (SLP, for short) (and hence the WLP).  From this, the following question arises naturally:

\begin{question}\label{qu:CI}
Does every complete intersection ideal $I\subseteq \mathbb{K}[x_1,\dots, x_n]$, in any number of variables, has the WLP (or even the SLP)?
\end{question}
This question should be considered as the main open problem concerning Lefschetz properties for complete intersections. In a step towards an answer to \Cref{qu:CI}, it seems reasonable to first try to extend the work in  \cite{HMNW} as follows.

\begin{question}
    Does every height three complete intersection has the SLP? Is this true for complete intersections of larger height, e.g., Artinian height four?
\end{question}

Even, more specific is the following problem, whose answer is also unknown.

\begin{question}
Let $I\subseteq \mathbb{K}[x_1,\dots, x_n]$ be a complete intersection generated by homogeneous polynomials of uniform degree $d\in \mathbb{N}$. Further assume that $\mathbb{K}[x_1,\dots, x_n]/I$ has the same monomial basis as $\mathbb{K}[x_1,\dots, x_n]/(x_1^d, x_2^d, \ldots, x_n^d)$. Does $\mathbb{K}[x_1,\dots, x_n]/I$  satisfy the SLP?
\end{question}

The last problem we are proposing for complete intersections is concerned with Lefschetz properties for non-standard graded
Artinian complete intersection and relates to an old conjecture by Almkvist \cite{Alm0,Alm1,Alm2} which we now explain. Let $p_{n,d}(t)$ be the generating function for the number of partitions of  integers smaller than or equal to $(d-1)n(n+1)/2$ with at most $n$ parts and each part repeated at most $d-1$ times.
It is known that
\[
p_{n,d}(t)=\prod_{i=1}^n\left(1+t^i+\cdots+t^{i(d-1)}\right).
\]

\begin{conjecture}[Almkvist '85]\label{conj:Almkvist}
 For every fixed $d\geq 2$, the polynomial $p_{n,d}(t)$ has unimodal coefficients for $n$ sufficiently large.
\end{conjecture}
Almkvist has already established his conjecture \cite{Alm2} for $3\leq d\leq 20$ (and also, oddly enough, for $d=100$ and $d=101$) using analytic techniques, including integrals and Tchebychev polynomials.

One possible approach towards \Cref{conj:Almkvist} is via Lefschetz properties. We now explain this in more detail.
Given integers $n,d\geq 1$, let $e_i=e_i(x_1,\ldots,x_n)$ be the \emph{$i^{th}$ elementary symmetric polynomial} in the variables $x_1,\ldots,x_n$, and let $\hat{e}_i=e_i(x_1^d,\ldots,x_n^d)$ be the $i^{th}$ elementary symmetric polynomials in the variables raised to the $d^{th}$ powers $x_1^d,\ldots,x_n^d$.  It is known that $\mathbb{K}[e_1,\ldots,e_n]$ is a polynomial ring, and that $\{\hat{e}_1,\ldots,\hat{e}_n\}\subseteq \mathbb{K}[e_1,\ldots,e_n]$ form a regular sequence. Consequently,  the quotient ring $A(n,d)=\frac{\mathbb{K}[e_1,\ldots,e_n]}{(\hat{e}_1,\ldots,\hat{e}_n)}$ is a (non-standard) graded Artinian complete intersection. It is easy to see that the Poincar\'e polynomial of $A(n,d)$ equals $p_{n,d}(t)$. Hence, an answer to the following question would solve \Cref{conj:Almkvist}.
\begin{question}
For which pairs $(n,d)$ does the ring $A(n,d)$ have the SLP?
\end{question}

{\bf Gorenstein algebras.}
The question of which Gorenstein algebras do possess the WLP (or even the SLP) is even more mysterious than the corresponding one for complete intersections. On one hand, as mentioned above,  all graded algebras in two or fewer variables have the WLP \cite[Proposition 4.4]{HMNW}. On the other hand, Ikeda and Boij showed that there \emph{exist} Artinian Gorenstein algebras of codimension $>3$ failing the WLP. In view of the fact, that complete intersections of codimension 3 do have the WLP, it is natural to ask the following question:
\begin{question}
Do all Artinian codimension three Gorenstein algebras have the WLP (or even the SLP)? If not, can one characterize those which do?
\end{question}
Building on work of Boij and Ikeda, one might also aim to characterize (or at least better understand) which Gorenstein algebras do have the WLP/SLP, e.g., in terms of properties of their Hilbert functions. However, it is not enough to only consider Gorenstein algebras with unimodal Hilbert functions (which is a consequence of the WLP), since, even in this more restricted setting, there exist examples of Gorenstein algebras that do not have the WLP. 

As for complete intersections, we now mention a more specific problem.

\begin{question}
Let $A = \bigoplus^c _{i=0} A_i $ be a standard graded Artinian Gorenstein algebra over a field of characteristic zero or greater than $c$ with embedding
dimension $n$. Suppose that the symmetric group $S_n$ acts on $A$ by permutation of the
variables.
Does $A$ have the SLP?
\end{question}
The above question has been answered in the affirmative for quadratic complete intersections.

\vskip 4mm
{\bf Jacobian ideals.}
Let $S=\mathbb C [x_0,\dots,x_n]$ be the polynomial ring in $n+1$ variables with coefficients in $\mathbb{C}$ and let $f \in S$ be a homogeneous polynomial of degree $d$ whose corresponding hypersurface $V=V(f)=\{x\in \PP^n~:~f(x)=0\}$ in $\PP^n$ is smooth. Let $J_f$ be the \emph{Jacobian ideal} of $f$, generated by the partial derivatives $f_i$ of $f$ with respect to $x_i$ for $i=0,\dots,n$. The graded algebra $M(f)= S/J_f$ is called the \emph{graded Milnor} (or \emph{Jacobian}) algebra of $f$. If $f$ is generic, it follows from  \cite{stanley}  and  \cite{watanabe}  that $M(f)$ has both the WLP and the SLP. Moreover, for any $f\in \mathbb C [x_0,\dots,x_3]$ of degree $d\le 6$ with $V(f)$ smooth we know that $M(f)$ satisfies the WLP (see \cite{GI} and \cite{BMMN}). This gives rise to the following question:

\begin{question}
Does $M(f)$ have the WLP (or the SLP) for any homogeneous form $f\in \mathbb C [x_0,\dots,x_n]$ with $V(f)\subset \PP^n$ smooth?
\end{question}

{\bf Monomial algebras.}
Though it is out of reach to provide a complete classification of all Artinian monomial algebras $A$ that have the WLP, one can still ask for a characterization for special classes.  Two natural subclasses come from looking at the end and at the beginning of a minimal free resolution of $A$.

An algebra $A$ is called \emph{level} if the  last free module in its minimal free resolution is concentrated in one degree. The rank of the last free module is called the {\em type} of $A$. It is natural to ask which monomial level algebras of a prescribed type have the WLP. It was shown in \cite{stanley} and \cite{watanabe} that for type 1 this is always true. Moreover,  in \cite{BMMNZ} it was proven that also monomial level algebras of type 2 in three variables always have the WLP, whereas for every other number of variables ($\geq 3$) and every possible type there exist monomial level algebras not having the WLP. This motivates the following question:
\begin{problem}
    Let $n$ and $t$ be positive integers. Find or bound the minimal degree $d$ such that there exists a monomial level algebra in $n$ variables of type $t$ whose last syszygy module is generated in degree $d$ and fails the WLP.
\end{problem}
We want to emphasize that obviously this question also makes sense for the SLP, for which even less is known in this setting.

Another class of monomial algebras, which is natural to consider are \emph{monomial almost complete intersections}. For these, several authors have studied the special case of three variables. It would be interesting to know what happens if we allow  one more generator. It is also not known what happens if the algebra is not required to be monomial but just an almost complete intersection.

Once more, the results that are known in characteristic zero change dramatically in characteristic $p$. A broad question is the following:
\begin{question}
 Let $I$ be a monomial ideal such that $R/I$ has the WLP in characteristic zero. What are the field characteristics in which $R/I$ fails to have the WLP?
\end{question}

{\bf Powers of linear forms.}  Many other natural algebras lend themselves to questions about the WLP or the SLP. A popular  instance is that the underlying ideal is generated by powers of linear forms (see, for instance, \cite{JMMW}). Most of the results achieved in this direction rely on a result of Emsalem and Iarrobino \cite{EI}, that translates the problem of whether the considered algebra has the WLP to one of studying sets of fat points in  a projective space.  In \cite{MMN}, the following conjecture was stated:

\begin{conjecture}\label{conj:MMN}
Let $R = \mathbb{K}[x_1,\dots,x_{2n+1}]$. Let $L \in R$ be a general linear form, and let $I = \langle x_1^d,\dots,x_{2n+1}^d,L^d \rangle$.
\begin{itemize}
\item[(i)] If $n=3$, the algebra $R/I$ fails the WLP if and only if  $d\geq 3$.
\item[(ii)] If $n\geq 4$, the algebra $R/I$ fails the WLP if and only if $d>1$.
\end{itemize}
\end{conjecture}

This conjecture has been solved for $d\ne 3$ and $n=3$  in \cite{MMN} and for all other cases in two steps by work of
Nagel and Trok \cite{NT}; and
Boij and Lundqvist \cite{BL}. So, we are let to pose the following problem:

\begin{problem}
    Let $R = \mathbb{K}[x_1,\dots,x_{m}]$. Let $L_1,\dots ,L_r \in R$ be  general linear forms, and let $I = \langle x_1^d,\dots,x_{m}^d,L_1^d ,\dots ,L_r\rangle$. Determine for which values of $m$, $d$ and $r$ the algebra $R/I$ fails t.he WLP
\end{problem}

\vskip 4mm
 {\bf  Gotzmann ideals.} A square-free monomial ideal $I\subseteq S=\mathbb{K} [x_0,\cdots x_n]$ is called
\emph{Gotzmann ideal} if and only if $$I=m_1(x_i~:~i\in J_1) +
m_1m_2(x_i~:~i\in J_2) + \cdots + m_1m_2\cdots m_s(x_i~:~i\in J_s)$$ for some square-free monomials
$m_1,\ldots,m_s$ and pairwise disjoint subsets $J_1,\dots,J_s$ of $\{0,1,\ldots,n\}$.

The above definition was introduced by Hoefel and Mermin in \cite{HM}
Herzog and Hibi \cite{HHCWL} showed that all Gotzmann ideals are componentwise linear, and  Bigdeli and Faridi~\cite{BF} established a connection between square-free Gotzmann and
Stanley-Reisner ideals of \emph{chordal complexes}~--~a large class of
componentwise linear ideals which can be defined via simplicial collapses. Open questions that arise from this are the following:
\begin{question}
 Which square free Gotzmann ideals satisfy the WLP?
 Do (some) Artinian reductions of Stanley-Reisner ideals of chordal complexes satisfy the WLP?
\end{question}

See \cite{TH} for some recent contributions to this question and related problems.

\section{Geometric aspects  related to Lefschetz properties}

In the last decade the failure or presence of the WLP has been connected to a large number
of geometric problems, that appear to be unrelated at first glance. For example, in \cite{MMO}, Mezzetti, Mir\'o-
Roig and Ottaviani proved that the failure of the WLP is related to the existence of projective varieties
satisfying at least one Laplace equation of order greater than 2.
Another connection with classical algebraic geometry has been discovered by
Maeno and Watanabe, who proved that an Artinian Gorenstein algebra $A=\mathbb{K}[x_1,\dots ,x_n]/\mathrm{Ann}(F)$ fails the SLP if one of the higher Hessians of $F$ is identically
zero \cite{MW}. This motivates the study of projective hypersurfaces with vanishing
Hessian or higher Hessians, a classical problem that goes back to Gordan and Noether.
In this section we will exhibit other geometric properties closely related to the failure or presence of the WLP.
As in the previous section we divide this section into several subsections.

\vskip 4mm
{\bf Reduced sets of points.} It is known that not every Artinian ideal is an Artinian reduction of the ideal of a reduced set of points, and certainly it is not a {\em general} Artinian reduction of the ideal of a reduced set of points.
An open question is the following:
\begin{question}\label{qu:GeneralReduction}
 Does a general Artinian reduction of a reduced, arithmetically Gorenstein set of points in $\PP^n$ has the WLP (provided that the  characteristic of the underlying field equals zero)?
\end{question}
A positive answer to this question would imply the algebraic $g$-conjecture  and yields a complete classification of the Hilbert functions of such sets of points. We want to emphasize that it is important to consider \emph{general} Artinian reductions since  examples of reduced arithmetically Gorenstein sets of points exist for which  a {\em special} Artinian reduction fails the WLP. A similar but less ambitious question whose positive answer would be implied by \Cref{qu:GeneralReduction} is the following:
\begin{question}
 Is the $h$-vector of  every reduced, arithmetically Gorenstein set of points unimodal (provided that the characteristic of the underlying field equals zero)?
\end{question}

 {\bf Artinian reductions of arithmetically Cohen-Macaulay curves.}
A question which is of the same flavour and completely wide open is the following:
\begin{question}
Do all irreducible arithmetically Cohen-Macaulay curves in $\mathbb P^4$ have an Artinian reduction that has the SLP?
\end{question}

{\bf Unexpected curves and complex line arrangements.}
Unexpected sets of points relate directly to failures of the SLP.
In particular, a finite set $Z\subseteq {\mathbb P}^n$ of $r$ points $p_1,\dots,p_r$ dual to hyperplanes
$H_1,\dots,H_r$ has an \emph{unexpected hypersurface cone} of degree $d$ with a general point
of multiplicity $d$ if and only if ${\mathbb C}[x_0,\ldots,x_n]/(L_1^d,\ldots,L_r^d)$
fails SLP in degree $d-1$ with range 1 (see \cite[Proposition 2.17]{HMNT}).

Unexpectedness in the plane often comes from line arrangements. In particular, there are four known kinds of arrangements of distinct lines $L_1,\ldots,L_r$ in ${\mathbb P}^2$ (over the complex numbers) such that whenever $L_i$ and $L_j$ meet there is a third line $L_k$ meeting at the same point. They are:
\begin{itemize}
\item 3 or more concurrent lines,
\item the Fermat  arrangements ($(x_0^r-x_1^r)(x_0^r-x_2^r)(x_1^r-x_2^r)$),
\item the Klein arrangement of 21 lines,
(this has 21 points where exactly 4 lines meet, 28 points where exactly 3 lines meet, and no other points where 2 or more lines meet),
\item the Wiman arrangement of 45 lines
(this has 36 points where exactly 5 lines meet, 45 points where exactly 4 lines meet, 120 points where exactly 3 lines meet and no other points where 2 or more lines meet).
\end{itemize}
For $r\geq 5$, the points dual to the lines of the Fermat arrangements have unexpected curves.
The points dual to the Klein and Wiman lines also have unexpected  curves (see \cite{CHMN}).
This naturally raises the following open problem:

\begin{question}
Are there other complex line arrangements having no points where exactly two lines meet?
\end{question}

{\bf Unexpected hypersurfaces and companion varieties.}
The existence of an unexpected curve was first established for the points of the $B_3$ root system. In \cite{Szp19} it was shown that the unexpectedness is directly related to the existence of Togliatti-type surfaces (having defective osculating spaces). Moreover, the author showed, in the case of the surface associated to the $B_3$ root system, that there exists another surface, called a \emph{companion surface}, which exhibits a number of interesting geometrical properties. This direction of study has been pursued in \cite{DGIMRSS} in case of Togliatti-type varieties and their companions associated to the $B_4$, $F_4$ and $H_3$ root systems and to certain Fermat-type configurations of points. This motivates the following problem:
\begin{problem}
Study Togliatti-type varieties and their companion surfaces for other root systems and find a Togliatti-type construction associated to configurations of points allowing multiple general fat points as described in \cite{Szp22}.
\end{problem}

{\bf Sets of points projecting to complete intersections}
A finite set $Z\subseteq {\mathbb P}^n$ is said to be \emph{geproci} if its image $Z_P$ under projection
from a general point $P$ to a fixed hyperplane $H$ is a complete intersection. An example
is a set of points $Z$ contained in a hyperplane $L$ which is already a complete intersection in $L$.
The only other examples of geproci sets known are in ${\mathbb P}^3$. In this case,
we say $Z$ is \emph{$(a,b)$-geproci} if $Z_p$ is the transverse intersection of a curve $C_{a,P}$ of degree $a$ with a curve $C_{b,P}$ of degree $b$ where $a\leq b$ (hence $|Z|=ab$).
Non-degenerate examples of such sets are given by \emph{$(a,b)$-grids}.
A set $Z$ of $ab$ points in ${\mathbb P}^3$ is called
an \emph{$(a,b)$-grid} if there is a set $A$ of $a$ skew lines and a set $B$ of $b\geq a$ skew lines,
such that each line in $A$ meets each line in $B$ in exactly one point (see \cite{POLITUS}).
By a theorem of \cite{CM}, if $Z$ is an $(a,b)$-grid with $3\leq a\leq b$, then
$C_{a,P}$ and $C_{b,P}$ are unexpected curves
and thus give failures of SLP (see the previous paragraph).

We define two more types of geproci sets. A \emph{halfgrid} is a geproci set $Z$
such that either $C_{a,P}$ or $C_{b,P}$ (but not both) is a union of lines for a general $P$.
A \emph{non-halfgrid} is a geproci set $Z$ such that neither $C_{a,P}$ nor $C_{b,P}$
is a union of lines for a general $P$.

Since degenerate geproci sets and grids are easy to construct and well-understood,
they are referred to as \emph{trivial} in the folowling. The obvious question is whether other geproci sets exist. More precisely, the following problem is open:
\begin{problem}
Find a non-trivial $(a,b)$-geproci set, where $C_{a,P}$ and $C_{b,P}$ is not unexpected, or show that no such example exists.
\end{problem}
 So far, not a single example of a non-trivial $(a,b)$-geproci set where $C_{a,P}$ and $C_{b,P}$ are not unexpected, is known.
Moreover, since every nontrivial geproci set known has at least one subset of 3 collinear points, the following problem is natural:
\begin{problem}
Prove or disprove (by giving a counterexample) that every non-trivial $(a,b)$-geproci set has a subset of $3$ collinear points.
\end{problem}
As a more general version, one can even consider the following problem:
\begin{problem}
Find an example of a linearly general geproci set, or show that no such example exists.
\end{problem}
In fact, all non-trivial geproci sets known have a rich struccture of linearly dependent subsets. Associating a matroid to these subsets, gives rise to Terao type problems, for instance:
\begin{question}
   Are geproci sets with isomorphic matroids projectively equivalent?  Is \emph{being geproci} a property of the associated matroids, i.e., is a set of points whose matroid is isomorphic to the matroid of a geproci set, geproci itself?
\end{question}


\section{Lefschetz properties in combinatorics}

 Lefschetz properties have very nice applications to ombinatorics. The main goal of this section is to gather several open problems concerning Lefschetz problems for simplicial complexes and more precisely for  their Stanley--Reisner ideals. As for the other sections, we divide this section into several topic areas.

We first fix some standard notation.
Let $S=\mathbb{K}[x_1,\dots,x_n]$ be the standard graded polynomial ring over a field $\mathbb{K}$ of characteristic zero.
Usually $\Delta $ will denote a simplicial complex on vertex set $[n]=\{1,\ldots,n\}$ and $I_\Delta$, respectively $\mathbb{K}[\Delta]$ will denote the corresponding Stanley--Reisner ideal respectively Stanley-Reisner ring. These, themselves, are not Artinian rings, but knowing Lefschetz properties of Artinian reductions of Stanley--Reisner rings is an important problem from both algebraic and combinatorial viewpoints.
One of the biggest progress on this topic is the following recent result of Adiprasito \cite{Adi} (see also \cite{APP,KaruX,PP})

\begin{thm}
If $\Delta$ is a simplicial $(d-1)$-sphere, then for a generic linear system of parameters  $\Theta$ of $S/I_\Delta$, the Artinian algebra $S/(I_\Delta+(\Theta))$ has the SLP.
\end{thm}

The special case of this theorem that $\Delta$ is the boundary complex of a simplicial polytope is a classical result due to Stanley and follows from the Hard Lefschetz theorem for toric varieties (see \cite[III \S 1]{Sta}).
It should also be mentioned that the theorem has a generalization to simplicial manifolds and doubly Cohen--Macaulay simplicial complexes (see \cite{APP}).
While these works give answers to many unsolved problems, still a variety of open problems on this topic remains.

{\bf Balanced spheres}
A simplicial complex on $[n]$ of dimension $(d-1)$ is said to be \emph{(completely) balanced} if there exists a function $c:[n] \to [d]$ such that $c(v)\ne c(u)$ for all edges $\{u,v\} \in \Delta$.
It should be noted that the function $c$ can be regarded as a proper \emph{coloring} of the $1$-skeleton, i.e., the graph of $\Delta$.
If $\Delta$ is balanced, then it is known that the sequence of linear forms $\theta_1,\dots,\theta_d$ defined by $\theta_i=\sum_{c(v)=i} x_v$ forms a linear system of parameters of $\Delta$, called a \emph{colored linear system of parameters} of $\Delta$ (see \cite[II \S4]{Sta}).
The following question is asked in \cite[Conjecture 1.3]{CJMN}.

\begin{problem}
Let $\Delta$ be a balanced simplicial $(d-1)$-sphere and $\Theta=\theta_1,\dots,\theta_d$ be a colored linear system of parameters of $\Delta$.
Does the algebra $S/(I_\Delta+(\Theta))$ have the SLP (or WLP)?
\end{problem}

{\bf Centrally symmetric spheres}
A simplicial complex $\Delta$ on vertex set $[n]$ is said to be \emph{centrally symmetric}
 (\emph{cs}, for short) if it admits a free $\mathbb Z/2\mathbb Z$-action, that is, there is a fixed-point free involution $\alpha:[n] \to [n]$ with $\Delta=\{\alpha(F)\mid F \in \Delta\}$.
For such a simplicial complex, one can choose a linear system of parameters $\Theta=\theta_1,\dots,\theta_d$ from $S^-=\{f \in S \mid \alpha(f)=-f\}$, where $\alpha(x_v)$ is defined to be $x_{\alpha (v)}$.

\begin{problem}
Let $\Delta$ be a centrally symmetric simplicial sphere and let $\Theta=\theta_1,\dots,\theta_d$ be a linear system of parameters for $\mathbb{K}[\Delta]$ that has been generically taken from $S^-$.
Does the algebra $S/(I_\Delta+(\Theta))$ have the SLP or WLP with respect to the linear form $w=x_1+ \cdots +x_n$?
\end{problem}

The answer is known to be yes for cs simplicial polytopes (see \cite[III \S8]{Sta} for more background).

\vskip 4mm
{\bf Flag spheres}
A simplicial complex $\Delta$ is said to be \emph{flag} if $I_\Delta$ has no generators of degree $\geq 3$. $h$-vectors of flag simplicial spheres have been of great interest in combinatorics and there are two famous conjectures on this topic, known as \emph{Charney-Davis Conjecture} and \emph{Gal's gamma non-negativity Conjecture}.
Motivated by these conjectures Hailun Zheng suggested the following question.

\begin{problem}
Let $\Delta$ be a flag simplicial $(d-1)$-sphere with $d \geq 5$,
let $\Theta=\theta_1,\dots,\theta_d$ be a generic linear system of parameters and let $w,w'$ be additional generic linear forms.
Is it true that the multiplication $\times w' : (S/(I_\Delta+(\Theta,w))_1 \to (S/(I_\Delta+(\Theta,w))_2$ is injective?
\end{problem}

{\bf Almost Lefschetz properties}
A special instance of a flag simplicial complex is the barycentric subdivisions $\mathrm{sd}(\Gamma)$ of a CW complex $\Gamma$.
If $\Gamma$ is a $(d-1)$-dimensional Cohen--Macaulay complex, then the following almost Lefschetz property is known (see \cite[Theorem 6.2]{MuraiYanagawa})
\begin{itemize}
\item[(i)] $\times w^{d-2k-1}: \big(S/(I_{\mathrm{sd}(\Gamma)}+(\Theta)) \big)_k \to  \big( S/(I_{\mathrm{sd}(\Gamma)}+(\Theta)) \big)_{d-1-k}$ is injective for $k \leq (d-1)/2$;
\item[(ii)] $\times w^{d-2k-1}: \big(S/(I_{\mathrm{sd}(\Gamma)}+(\Theta)) \big) _{k+1} \to \big( S/(I_{\mathrm{sd}(\Gamma)}+(\Theta)) \big)_{d-k}$ is surjective for $k \leq (d-1)/2$.
\end{itemize}
Here, $\Theta$ is a generic linear system of parameters, and $w$ is a generic linear form.
This result in particular implies that $S/(I_{\mathrm{sd}(\Gamma)}+(\Theta))$ has the WLP if $d$ is odd.

\begin{question}
Is there a Cohen--Macaulay CW complex $\Gamma$ such that $S/(I_{\mathrm{sd}(\Gamma)}+(\Theta))$ does not have the weak Lefschetz property for any linear system of parameters $\Theta$?
\end{question}

As a related vague question, one could also ask
\begin{question}
Are there other classes of simplicial complexes (different from barycentric subdivisions) that satisfy (i) and (ii) above?
\end{question}
While property (i) is known to hold  for more general subdivisions \cite[Theorem 46]{AG}, property (ii) seems a bit mysterious at this moment.

\vskip 4mm
 {\bf Nagata idealizations.}
 The idealization construction by Nagata has been used in many cases to
construct Artinian Gorenstein algebras failing the WLP, for example
the famous example by R. Stanley. In terms of Macaulay inverse systems,
we can construct the idealization of the canonical module of an
Artinian level algebra with inverse system $\langle
F_1,F_2,\dots,F_s\rangle$ as the Gorenstein algebra with inverse
system $\langle \sum x_iF_i\rangle$ where $x_1,x_2,\dots,x_s$ are new
variables. Here the forms $F_1,\dots, F_s$ are assumed to be
homogeneous of degree $d$ in variables $y_1,y_2,\dots,y_n$.
Preliminary studies of the Gorenstein algebra given by
$\langle \sum x_i^eF_i\rangle$ show that for $e\ge d$ it satisfies the WLP for any forms $F_i$ while for $e<d$ we have examples failing the WLP and examples satisfying the WLP. THerefore  we are in front of an interesting problem:

\begin{problem}
  To determine whether the Gorenstein algebra given by $\sum x_i^eF_i$ has the WLP and/or the SLP.
\end{problem}

As a first contribution to this last problem the reader can look at \cite{AADFIMMMRN} and \cite{FMMR}.

\section{Lefschetz properties and Jordan types of Artinian algebras}

Let $\mathbb{K}$ be any field, and let $A$ be an  Artinian algebra, quotient of the polynomial ring $\mathbb{K}[x_1,\ldots,x_n]$ or of the local regular ring $\mathbb{K}\{x_1,\ldots,x_n\}$. Let ${\mathfrak{m}=(x_1,\ldots,x_n)}$ be the maximal ideal of $A$. Let $M$ be a finite module over $A$, and let ${\ell\in\mathfrak{m}}$. The \emph{Jordan type} ${P_{\ell,M}}$ of $\ell$ is the partition of $\dim_\mathbb{K}M$ giving the sizes of the Jordan blocks in a Jordan canonical form for the linear map ${\times\ell:M\to M}$, defined by the multiplication by $\ell$. There are several refinements of this invariant. We refer to \cite{IMM23} in this volume for their definitions and a discussion of their properties.  In particular, when $M$ is a graded module over an Artinian algebra $A$ we define the generic linear Jordan type $P_{1,M}$ as the Jordan type $P_{\ell,M}$ for a generic linear form $\ell\in A_1$.

The Jordan type has two main strengths. First, in the graded case, it is a finer invariant than the Lefschetz properties: we can determine if a linear form $\ell$ satisfies the WLP or the SLP on $A$ from its Jordan type and the Hilbert function $H=H(A)$: that is, $\ell$ has the SLP if its Jordan type is the conjugate $H^\vee$; and $\ell$ has the WLP if the number of parts of the Jordan type is the maximum value of $H$.  Second, the Jordan type is well-defined also on non-{graded} algebras, allowing us to generalize the definitions of WLP and SLP to the non-{graded} case. As in the previous sections we distinguish several topics:

\vskip 4mm
{\bf Linear Jordan type and contiguous Jordan type.}	Let $A$ be standard graded, if $H=H(A)$ is unimodal (and has a single maximum value with no dips) and if the generic Jordan type is the conjugate partition of $H$ (i.e., $A$ is strong Lefschetz) then the generic linear Jordan type is equal to $H^\vee$, and therefore they are the same \cite[Proposition 2.14]{IMMM22a}. \par
\begin{question} Is this equality still true when $H(A)$ is not unimodal, or when the generic Jordan type is not the conjugate partition of the Hilbert function? More generally, under which conditions on a graded module $M$ over a graded Artinian algebra $A$ do their generic  and generic linear Jordan types satisfy ${P_M=P_{1,M}}$? (See Questions 1.1 and 2.56 in \cite{IMMM22a}.)
\end{question}
We can decompose $M$ as the direct sum of ${\mathbb K}[\ell]$ modules~--~called $\ell$-strings~--~whose lengths are given by the Jordan partition ${P_{\ell,M}}$,
The Jordan degree type JDT $P_{\mathrm{deg},\ell,M}$ of a graded module $M$ adds to the Jordan type the information about the initial degree of the $\ell$-strings in such a decomposition of $M$ as  module over $\mathsf{k}[\ell]$, this depends only on the pair $(\ell,M)$. However, a problem is that this definition of JDT does not generalize to non-graded Artinian algebras or modules \cite{IMS}.\par
When the Hilbert function $H$ of a graded Artin module is non-unimodal, we can define a contiguous-Jordan type $P_{c,\ell}(H)$, and a contiguous-Jordan degree type
$P_{c,\deg,\ell}(H)$: using the bar graph of $H(M)$ \cite[Definition 2.28]{IMMM22a}.
\begin{question} For which (non-unimodal) Hilbert functions $H$ occurring for graded Artinian algebras $A$, can we find pairs $(A,\ell)$ with $H(A)=H$ and the Jordan type and Jordan degree types of $A$ agree with the contiguous Jordan or Jordan-degree type of $H$?
\end{question}

\vskip 4mm
{\bf Jordan type of the initial ideal.}
Recall that the {\em initial ideal} $\mathrm{init}(I)$ of an ideal $I\subseteq R=\mathbb{K}[x_1,\ldots,x_n]$ is the ideal generated by the leading terms of elements of $I$ under a fixed monomial order.
From a result in  \cite{C} and \cite{W} we know that if $R/(\mathrm{in} (I) )$ has the WLP (resp. the SLP), then the same holds for $R/I$. More precisely,  for a generic linear form $\ell$ we have that for all $j,k$ the following inequality is true
$$\dim_\mathbb{K}(R/(I,\ell^k)_j\le \dim_\mathbb{K}(R/(\mathrm{in} I,\ell^k))_j.$$
\begin{question}
	How does the generic Jordan degree type of a standard {graded} Artinian algebra $A=R/I$ behave under projection to the quotient $R/\mathrm{in}( I)$? \par
	Is there a pair $(A,\ell)$, where $A$ is a standard {graded} Artinian algebra and ${\ell\in A_1}$, such that the Jordan degree type $P_{\deg,\ell,A}$ cannot occur for a pair $(A',\ell')$, if $A'$ is defined by a monomial ideal?
\end{question}

\vskip 4mm
{\bf Jordan type for non-graded algebras.}
	Lately, some refinements of Jordan type have been introduced for non-graded Artinian algebras: namely sequential Jordan type, Loewy sequential Jordan type, or double  sequential Jordan type \cite{IMS}. Each of these is semicontinuous in a family of algebras having fixed Hilbert function.  The semicontinuity of Jordan type has been used to show that certain families $\mathrm{Gor}(H)$ of Artinian Gorenstein algebras with given Hilbert function have several irreducible components, beginning in codimension three \cite{IMM}.
	\begin{question} Determine Hilbert function sequences $H$ for non-graded Artinian algebras such that the family $Z(H)$ of all algebras with Hilbert function $H$ has several irreducible components  $\Xi_1,\Xi_2$- whose general elements have the same generic Jordan type, but which differ in one of the refined invariants?
\end{question}
\begin{question}
	Let $\,{0\to L\to M\to N\to0}$\, be an exact sequence of finite modules over $A$. How can we compare the generic Jordan types $P_L$, $P_M$, and $P_N$? Under what conditions could we have additivity ${P_M = P_L + P_N }$, in a suitable sense for the Jordan types? The same question can be extended to Jordan degree type, in the graded case, or to sequential Jordan type, Loewy sequential Jordan type, or double  sequential Jordan type, in the non {graded} case.
\end{question}

The next question has to do with the symmetric decomposition of the associated graded algebra of a non-graded Gorenstein Artinian algebra (see \cite{IMM}). We know, by an example of Chris McDaniel \cite{IMS} that it is not always possible to find a Jordan basis $B$ for the multiplication by an element ${\ell\in\mathfrak{m}}$ that agrees with the Hilbert function of $A$, in the sense that ${B\cap\mathfrak{m}^i}$ is a basis for the vector space ${\mathfrak{m}^i}$ for every $i$. We also know that there is always a pre {Jordan} basis that agrees with the Hilbert function in this sense \cite{IMS}. So we are let to pose the following question:

\begin{question}
	 If the Artinian algebra $A$ is Gorenstein, can we always find a pre {Jordan} basis $B$ that agrees with the $Q(a)$ decomposition, in the sense that ${B\cap\mathfrak{m}^i\cap(0:\mathfrak{m}^b)}$ is a basis for the vector space ${\mathfrak{m}^i\cap(0:\mathfrak{m}^b)}$ for every $i$ and $b$?
\end{question}

\vskip 4mm
{\bf Jordan type of graded Artinian Gorenstein algebras.}
Graded Artinian Gorenstein (AG) algebras satisfy Poincar\'e pairing. This allows us to determine their Lefschetz properties  by the ranks of fewer multiplication maps. For instance,  a graded AG algebra $A$ satisfies the WLP if the multiplication map by a generic linear form has maximal rank in degree  $\lfloor\frac{d}{2}\rfloor$, where $d$ is the socle degree of $A$ (see \cite{MMU} where it is also shown that WLP can be sensitive to characteristic).  In codimension three, when the characteristic of $\mathbb{K}$ is zero, each standard-graded complete intersection algebra has the WLP \cite{HMNW}; some have conjectured that all codimension three standard graded AG algebras have the WLP (see \cite{BMMNZ} where it is proven that this conjecture depends on showing the special cases of compressed algebras).  It is open which graded AG algebras of codimension three have the SLP \cite{AAISY}. There are families of graded AG algebras of codimension four that fail the WLP \cite[Theorem 1.6]{AbS}. There are AG algebras in codimension five whose Hilbert functions are non-unimodal, and thus are also not WLP. The question of whether there are complete intersection algebras in codimension four (or even higher) that fail WLP has been asked by many researchers and is open.

\begin{question}
Determine the number of parts of generic linear Jordan types for the families of graded Artinian Gorenstein algebras failing the WLP.  More generally,  determine the generic linear Jordan type for a graded Artinian Gorenstein algebra that fails the SLP.
\end{question}

\begin{question}
Given a sequence $H$ that can occur as the Hilbert function of a graded Artinian Gorenstein algebra, can we find potential Jordan types or Jordan degree types? (Known in codimenson two \cite{AIK}).
\end{question}
\begin{question}
In \cite{AIK},  all possible linear Jordan types of complete intersection algebras of codimension two having a fixed Hilbert function is listed.  For every such Hilbert function $H$,  can we determine the pairs of Jordan types $P$ and $Q$ that occur for $(\ell_1,A)$ and $(\ell_2,A)$ where $H(A)=H$?
\end{question}


\section*{Acknowledgements}
Many people contributed to this chapter by supplying problems and we are thankful to all of them. Special thanks go to Nasrin Altafi, Brian Harbourne, Anthony Iarrobino, Pedro Marcua Marques, Satoshi Murai, Tomasz Szemberg and Justyna Szpond.


\end{document}